\documentclass[a4paper, 10pt,oneside]{amsart}
\usepackage{amsmath, amssymb, amsfonts}
\usepackage[all]{xy}
\newtheorem{tm}{Theorem}
\newtheorem{lm}[tm]{Lemma}
\newtheorem{prop}[tm]{Proposition}

\newenvironment{dokaz}
{\noindent\emph{Proof:}\ }
{\hfill $\blacksquare$}

\newcommand{\Z}
{{\mathbb Z}}
\newcommand{\N}
{{\mathbb N}}
\newcommand{\C}
{{\mathbb C}}

\newcommand{\g}
{{\mathfrak g}}

\newcommand{\gk}
{\hat{{\mathfrak g}}}
\newcommand{\gt}
{\tilde{{\mathfrak g}}}
\newcommand{\hk}
{\hat{{\mathfrak h}}}
\newcommand{\hz}
{\hat{\mathfrak h}_\Z}
\newcommand{\he}
{{\mathfrak h}^e}
\newcommand{\nt}
{\tilde{{\mathfrak n}}}
\newcommand{\h}
{{\mathfrak h}}
\newcommand{\n}
{{\mathfrak n}}
\newcommand{\gsl}
{{\mathfrak sl}}
\newcommand{\nop}
{\substack{\circ \\ \circ}}
\newcommand{\End}
{\textrm{End}\,}

\begin{document}

\author{Goran Trup\v{c}evi\'{c}}
\title{Combinatorial bases of Feigin-Stoyanovsky's type subspaces of level $1$ standard modules for $\tilde{\gsl}(\ell+1,\C)$} 
\address{Department of Mathematics, University of Zagreb, Bijeni\v cka 30, Zagreb, Croatia}
\curraddr{}
\email{gtrup@math.hr}
\thanks{}
\subjclass[2000]{Primary 17B67; Secondary 17B69, 05A19.\\ \indent Partially supported by the Ministry of Science and Technology of the Republic of Croatia, Project ID 037-0372794-2806}
\date{}
\dedicatory{}

\begin{abstract}

Let $\tilde{\mathfrak g}$ be an affine Lie algebra of type
$A_\ell^{(1)}$. Suppose we're given a $\mathbb Z$-gradation of
the corresponding simple finite-dimensional Lie algebra ${\mathfrak
g}={\mathfrak g}_{-1}\oplus{\mathfrak g}_0 \oplus {\mathfrak
g}_1$; then we also have the induced $\mathbb Z$-gradation of
the affine Lie algebra
$$\tilde{\mathfrak g}=\tilde{\mathfrak g}_{-1} \oplus
\tilde{\mathfrak g}_0 \oplus \tilde{\mathfrak g}_1.$$ 

Let
$L(\Lambda)$ be a standard module of level $1$. 
Feigin-Stoyanovsky's type subspace $W(\Lambda)$ is the $\tilde{\mathfrak
g}_1$-submodule of $L(\Lambda)$ generated by the highest-weight
vector $v_\Lambda$,
$$W(\Lambda)=U(\tilde{\mathfrak g}_1)\cdot v_\Lambda\subset L(\Lambda).$$
We find a combinatorial basis of $W(\Lambda)$ given
in terms of difference and initial conditions. Linear
independence of the generating set is proved  inductively by using
coefficients of intertwining operators. A basis of $L(\Lambda)$ is obtained as an ``inductive limit'' of the basis of $W(\Lambda)$. 
\end{abstract}

\maketitle


\section{Introduction}

Let $\g$ be a simple complex Lie algebra, $\h\subset\g$ its Cartan
subalgebra, $R$ the corresponding root system. Then one has a root
decomposition $\g=\h+\sum_{\alpha\in R}\g_\alpha$. Fix root vectors
$x_\alpha\in\g_\alpha$. Let
\begin{equation}\label{ZGradG_jed}\g=\g_{-1}\oplus \g_0 \oplus
\g_1\end{equation} be a $\Z$-gradation of $\g$, where
$\h\subset\g_0$. All such gradations are obtained by choosing some
minuscule coweight $\omega\in\h$. Denote by $\Gamma\subset R$ a set
of roots such that $\g_1=\sum_{\alpha\in\Gamma}\g_\alpha=\sum_{\omega(\alpha)=1}\g_\alpha$.

 Affine
Lie algebra associated with $\g$ is $\gt=\g\otimes
\C[t,t^{-1}]\oplus \C c \oplus \C d$, where $c$ is the canonical
central element, and $d$ the degree operator. Elements
$x_\alpha(n)=x_\alpha\otimes t^n$ are fixed real root vectors.
Gradation of $\g$ induces analogous $\Z$-gradation of $\gt$:
$$\gt=\gt_{-1}\oplus \gt_0 \oplus \gt_1,$$
where $\gt_1=\g_1\otimes\C[t,t^{-1}]$ is a commutative Lie
subalgebra with a basis
$$\{x_\gamma(j)\,|\,j\in\Z,\gamma\in\Gamma\}.$$

Let $L(\Lambda)$ be a standard $\gt$-module of level $k=\Lambda(c)$,
with a fixed highest weight vector $v_\Lambda$. A Feigin-Stojanovsky's
type subspace  is a $\gt_1$-submodule of $L(\Lambda)$ generated with
$v_\Lambda$,
$$W(\Lambda)=U(\gt_1)\cdot v_\Lambda\subset L(\Lambda).$$
This is similar to the notion of principal subspace introduced in
[FS] where, instead of $\Z$-gradation \eqref{ZGradG_jed}, one
considers triangular decomposition of $\g$ and from it derived
decomposition of $\gt$; in the case $\g=\gsl(2,\C)$, these two
definitions are equivalent.

We would like to find a {\em monomial basis} of $W(\Lambda)$,
i.e. a basis consisting of vectors $x(\pi)v_\Lambda$, where $x(\pi)$
are monomials in basis elements
$\{x_\gamma(-j)\,|\,j\in\N,\gamma\in\Gamma\}$. 

The problem of finding monomial bases is a part of
Lepowsky-Wilson's program to study representations of affine Lie
algebras by means of vertex-operators and to obtain
Rogers-Ramanujan-type combinatorial bases of these
representations ([LW], [LP], [MP]).

Principal subspaces of standard $\gk$-modules were
introduced in [FS]. These subspaces are generated by the
affinization of the nilpotent subalgebra $\n_+$ of $\g$ from the
triangular decomposition $\g=\n_- \oplus \h \oplus \n_+$. B.Feigin and
A.Stoyanovsky described the dual space of the principal subspace for
$\g=\gsl(2,\C)$ and $\gsl(3,\C)$ in terms of symmetric polynomial
forms satisfying certain conditions, and calculated its
character. In the $\gsl(2,\C)$-case, they also described the
dual in a geometric way, recovering in this way the Rogers-Ramanujan and Gordon
identities.

Principal subspaces were studied further by G.Georgiev in [G]. He
constructed combinatorial bases and calculated characters of
principal subspaces for certain representations of
$\gsl(\ell+1,\C)$. In the proof of linear independence, Georgiev
used intertwining operators from [DL].

Also by using intertwining operators, S.Capparelli, J.Lepowsky and A.Milas
 in [CLM1,2] obtained Rogers-Ramanujan and Rogers-Selberg recursions
 for characters of principal subspaces for $\gsl(2,\C)$. As a
 continuation of the program laid out in [CLM1,2], C.Calinescu
 obtained systems of recursions for characters of principal subspaces
of level $1$ standard modules for $\gsl(\ell+1,\C)$ ([C1]) and of
certain higher-level standard modules for $\gsl(3,\C)$ ([C2]). By solving
these recursions they also established formulas for characters of
these subspaces. Furthermore, in [CalLM1,2], Calinescu, Lepowsky and Milas
provided new proofs of presentations of principal subspaces for $\gsl(2,\C)$.

Feigin-Stoyanovsky's type subspace $W(\Lambda)$  was implicitly studied
in [P1] and [P2], where M.Primc constructed a combinatorial basis of
this subspace. By translating the basis of $W(\Lambda)$ by a certain
Weyl group element, and then taking a inductive limit, he obtained a basis of the whole $L(\Lambda)$.
This  was done in [P1] for $\g=\gsl(\ell+1,\C)$ and a particular
choice of gradation \eqref{ZGradG_jed}, and for any dominant integral weight $\Lambda$. For
any classical simple Lie algebra and any possible gradation
\eqref{ZGradG_jed}, combinatorial bases were constructed in [P2], but only for 
basic modules $L(\Lambda_0)$.

In the particular $\gsl(\ell+1,\C)$ case studied in [P1], the basis of
$W(\Lambda)$ is parameterized by combinatorial objects called
$(k,\ell+1)$-admissible configurations. These objects were introduced and further studied in [FJLMM] and [FJMMT], where different formulas for the character of $W(\Lambda)$ were obtained. 

The hardest part of constructing the combinatorial basis of
$W(\Lambda)$ is a proof of linear independence of a reduced
spanning set. This was proved in [P1] by using Schur functions,
while in [P2] this was proved by using the crystal base character
formula [KKMMNN]. In [P3], Primc used Capparelli-Lepowsky-Milas'
approach via intertwining operators and a description of the 
basis from [FJLMM] to give a simpler proof of linear independence of
the basis of $W(\Lambda)$ constructed in [P1]. It seems
that this should be the way to obtain a proof in other cases as
well.

In this paper we extend these results to any possible $\Z$-gradation
of $\g=\gsl(\ell+1,\C)$ and all level $1$ standard modules. In [T]
we will further extend this to standard modules of any higher level,
obtaining a combinatorial basis parameterized by a
certain generalization of $(k,\ell+1)$-admissible configurations.

Let $\delta=\{\alpha_1,\dots,\alpha_\ell\}$ be a basis of the root system $R$ for
$\g=\gsl(\ell+1,\C)$, and $\{\omega_1,\dots,\omega_\ell\}$ the
corresponding set of fundamental weights. We identify $\h$ and $\h^*$ in the usual way and fix a fundamental weight
$\omega=\omega_m$. Set $$\Gamma=\{\gamma\in
R\,|\,\langle\gamma,\omega\rangle=1\}=\{\gamma_{ij}\,|\, i=1,\ldots,m; j=m,\ldots,\ell\},$$
where $$ \gamma_{ij}=\alpha_i+\cdots+\alpha_m+\cdots+\alpha_j.$$ Set
 $${\mathfrak g}_{\pm1}  =
\sum_{\alpha \in \pm \Gamma}\, {\mathfrak g}_\alpha,\ {\mathfrak
g}_0 = {\mathfrak h} \oplus\sum_{\langle\alpha,\omega\rangle=0}\,
{\mathfrak g}_\alpha,$$ then $$ \mathfrak g  = \mathfrak g_{-1}
\oplus \mathfrak g_0 \oplus \mathfrak g_1$$ is a $\Z$-gradation of
$\g$. The set $\Gamma$ is called {\em the set of colors}. For $\gamma\in \Gamma$, we say that a fixed basis element $x_\gamma\in\g_\gamma$ {\em is of the color} $\gamma$. The
set of colors $\Gamma$ can be pictured as a rectangle with
row indices $1,\dots,m$ and column indices $m,\dots,\ell$ (see figure \ref{Gamma_fig}).

\begin{figure}[ht] \caption{} \label{Gamma_fig}
\begin{center}\begin{picture}(245,140)(-25,-10) \thicklines
\put(0,0){\line(1,0){180}} \put(0,0){\line(0,1){120}}
\put(180,120){\line(0,-1){120}} \put(180,120){\line(-1,0){180}}
\put(-25,117){$\Gamma$:} \put(-6,111){$\scriptstyle 1$}
\put(-6,99){$\scriptstyle 2$} \put(-8,3){$\scriptstyle m$}
\put(3,-8){$\scriptstyle m$} \put(12,-8){$\scriptstyle m
{\scriptscriptstyle +} 1$} \put(174,-8){$\scriptstyle \ell$}

 \linethickness{.075mm} \multiput(0,108)(4,0){45}{\line(1,0){2}}
\multiput(0,96)(4,0){45}{\line(1,0){2}}
\multiput(0,12)(4,0){45}{\line(1,0){2}}
\multiput(12,0)(0,4){30}{\line(0,1){2}}
\multiput(24,0)(0,4){30}{\line(0,1){2}}
\multiput(168,0)(0,4){30}{\line(0,1){2}}
\multiput(0,54)(4,0){45}{\line(1,0){2}}
\multiput(0,66)(4,0){45}{\line(1,0){2}}
\multiput(90,0)(0,4){30}{\line(0,1){2}}
\multiput(102,0)(0,4){30}{\line(0,1){2}}

\put(-6,57){$\scriptstyle i$} \put(93,-8){$\scriptstyle j$}
\put(90,57){$\gamma_{ij}$}
\end{picture}\end{center}
\end{figure}

Fix a fundamental weight $\Lambda_i, i=0,\dots,\ell$ of $\gt$. Let
$L(\Lambda_i)$ be the standard module with highest weight
$\Lambda_i$, and $v_i$ the highest weight vector of $L(\Lambda_i)$.

We find a basis of the Feigin-Stoyanovsky's type subspace
$W(\Lambda_i)$ consisting of {\em monomial vectors}
$$\{x_{\gamma_1}(-n_1) \cdots x_{\gamma_t}(-n_t) v_i\,|\,
t\in\Z_+;\gamma_j\in \Gamma, n_j\in\N\}$$ whose {\em monomial parts}
\begin{equation}
\label{MonomPart_def} x_{\gamma_1}(-n_1) \cdots x_{\gamma_t}(-n_t)
\end{equation}
satisfy certain combinatorial conditions, called {\em difference}
and {\em initial conditions}. By difference conditions, colors of
elements of degree $-j$ and $-j-1$ in a monomial
\eqref{MonomPart_def} lie on a diagonal path in $\Gamma$ as pictured
on the figure \ref{DCSeq_fig}.

\begin{figure}[ht] \caption{} \label{DCSeq_fig}
\begin{center}\begin{picture}(200,140)(-8,-10) \thicklines
\put(0,0){\line(1,0){180}} \put(0,0){\line(0,1){120}}
\put(180,120){\line(0,-1){120}} \put(180,120){\line(-1,0){180}}
\put(-6,111){$\scriptstyle 1$} \put(-6,99){$\scriptstyle 2$}
\put(-8,3){$\scriptstyle m$} \put(3,-8){$\scriptstyle m$}
\put(174,-8){$\scriptstyle \ell$}

\linethickness{.075mm} \multiput(88,72)(4,0){23}{\line(1,0){2}}
\multiput(88,72)(0,4){12}{\line(0,1){2}}

\thinlines \put(96,76){$\scriptscriptstyle \bullet$}
\put(98,78){\line(1,1){12}} \put(108,88){$\scriptscriptstyle
\bullet$} \put(132,94){$\scriptscriptstyle \bullet$}
\put(134,96){\line(1,1){12}} \put(144,106){$\scriptscriptstyle
\bullet$} \put(110,90){\line(4,1){24}} \put(108,97){$\scriptstyle
(-j)$}

\put(60,10){$\scriptscriptstyle \circ$}
\put(84,16){$\scriptscriptstyle \circ$}
\put(108,40){$\scriptscriptstyle \circ$}
\put(120,64){$\scriptscriptstyle \circ$}
\put(62,12){\line(4,1){24}} \put(86,18){\line(1,1){24}}
\put(110,42){\line(1,2){12}} \put(96,22){$\scriptstyle (-j-1)$}
\end{picture}\end{center}
\end{figure}

So, if a monomial \eqref{MonomPart_def} has elements of degrees $-j$
and $-j-1$ of colors $\gamma_{r_1 s_1},\dots,\gamma_{r_t s_t}$ and
$\gamma_{r_1' s_1'},\dots,\gamma_{r_{t'}' s_{t'}'}$, respectively,
then
$$r_1<r_2<\dots<r_t \textrm{\ and\ } s_1>s_2>\dots>s_t,$$
and, similarly,
$$r_1'<r_2'<\dots<r_{t'}' \textrm{\ and\ } s_1'>s_2'>\dots>s_{t'}'.$$
Also, $$r_t<r_1' \textrm{\ or\ } s_t>s_1'.$$

Initial conditions on monomials \eqref{MonomPart_def} require that
diagonal path of colors of elements of degree $-1$ lie below the
$i$-th row, in case $1\leq i\leq m$, or left of the $i$-th column,
in case of $m\leq i\leq \ell$, as it is pictured on the figure
\ref{ICSeq_fig}.

\begin{figure}[ht] \caption{} \label{ICSeq_fig}
\begin{center}
\begin{picture}(365,140)(-8,-10) \thicklines
\put(0,0){\line(1,0){160}} \put(0,0){\line(0,1){120}}
\put(160,120){\line(0,-1){120}} \put(160,120){\line(-1,0){160}}
\put(-6,111){$\scriptstyle 1$} \put(-6,99){$\scriptstyle 2$}
\put(-8,3){$\scriptstyle m$} \put(3,-8){$\scriptstyle m$}
\put(154,-8){$\scriptstyle \ell$}

\put(182,0){\line(1,0){160}} \put(182,0){\line(0,1){120}}
\put(342,120){\line(0,-1){120}} \put(342,120){\line(-1,0){160}}
\put(176,111){$\scriptstyle 1$} \put(176,99){$\scriptstyle 2$}
\put(174,3){$\scriptstyle m$} \put(185,-8){$\scriptstyle m$}
\put(336,-8){$\scriptstyle \ell$}

\linethickness{.075mm} \multiput(0,72)(4,0){40}{\line(1,0){2}}
\multiput(310,0)(0,4){30}{\line(0,1){2}} \put(-6,75){$\scriptstyle
i$} \put(313,-8){$\scriptstyle i$}

\multiput(2,72)(3,0){53}{\line(0,1){48}}
\multiput(310,3)(0,3){39}{\line(1,0){32}}

\thinlines \put(40,10){$\scriptscriptstyle \bullet$}
\put(42,12){\line(1,1){24}} \put(64,34){$\scriptscriptstyle
\bullet$} \put(66,36){\line(1,1){24}}
\put(88,58){$\scriptscriptstyle \bullet$}
 \put(62,26){$\scriptstyle (-1)$}

\thinlines \put(242,22){$\scriptscriptstyle \bullet$}
\put(244,24){\line(2,1){24}} \put(266,34){$\scriptscriptstyle
\bullet$} \put(268,36){\line(1,2){24}}
\put(290,82){$\scriptscriptstyle \bullet$}
 \put(264,26){$\scriptstyle (-1)$}

\end{picture}\end{center}
\end{figure}

Difference conditions on monomials are obtained by observing
relations between fields $x_\gamma(z),\gamma\in\Gamma$ on
$L(\Lambda_i)$, while initial conditions follow from the obvious
requirement that elements of degree $-1$ mustn't annihilate the
highest weight vector $v_i$.

By observing configurations of colors of elements of degrees $-1$
and $-2$, one is able to construct coefficients of suitable
intertwining operators between standard modules that would either send
basis elements of one module to basis elements of the other
module, or it would anihilate them. These operators are then used for the inductive proof of linear independence.

Thus we are able to prove the main result of this work

\vskip 1ex
\noindent
{\bf Theorem \ref{FSbaza_tm}}\quad {\em
Let $L(\Lambda_i)$ be a standard module of level $1$. Then the set of
monomial vectors $x_{\gamma_1}(-n_1) \cdots x_{\gamma_t}(-n_t) v_i$
whose monomial part satisfies difference and initial  conditions, is
a basis of $W(\Lambda_i)$.}

\section{Affine Lie algebras}

For $\ell\in \N$, let $$\g=\gsl(\ell+1,\C),$$ a simple Lie algebra
of the type $A_\ell$. Let $\h\subset\g$ be a Cartan subalgebra of
$\g$ and $R$ the corresponding root system. Fix a basis
$\Pi=\{\alpha_1,\dots,\alpha_\ell\}$ of $R$. Then we have the
triangular decomposition $\g=\n_-\oplus \h \oplus \n_+$. By $R_+$
and $R_-$ we denote sets of positive and negative roots, and let
$\theta$ be the maximal root. Let $\langle
x,y\rangle=\textrm{tr\,}xy$ be a normalized invariant bilinear form on
$\g$; via $\langle \cdot,\cdot\rangle$ we have an identification
$\nu:\h\to\h^*$. For each root $\alpha$ fix a root vector $x_\alpha\in\g_\alpha$.

Let $\{\omega_1,\dots,\omega_\ell\}$ be the set of
fundamental weights of $\g$,
$\langle\omega_i,\alpha_j\rangle=\delta_{ij},\,i,j=1,\dots,\ell$. Denote by $Q=\sum_{i=1}^\ell\Z \alpha_i$
the root lattice, and by $P=\sum_{i=1}^\ell\Z \omega_i$ the weight lattice
of $\g$. 

Denote by $\gt$ the associated affine Lie algebra
$$\gt=\g\otimes \C[t,t^{-1}]\oplus \C c \oplus \C d.$$
Set $x(j)=x\otimes t^j$ for $x\in\g,j\in\Z$. Commutation relations
are then given by
\begin{eqnarray*}
\ [c,\gt] & = & 0, \\
\ [d,x(j)] & = & j x(j), \\
\ [x(i),y(j)] & = & [x,y](i+j)+ i\langle x,y\rangle \delta_{i+j,0}c.
\end{eqnarray*}


Set $\he=\h\oplus\C c \oplus\C d,\, \nt_\pm=\g\otimes t^{\pm 1}\C
[t^{\pm 1}]\oplus \n_\pm$. Then we also have the triangular
decomposition $\gt=\nt_-\oplus\he\oplus\nt_+$.

Let $\hat{\Pi}=\{\alpha_0,\alpha_1,\dots,\alpha_\ell\}\subset
(\he)^*$ be the set of simple roots of $\gt$. Usual extensions of bilinear
forms $\langle\cdot,\cdot\rangle$ onto $\he$ and $(\he)^*$ are
denoted by the same symbols (we take $\langle c,d \rangle=1$).
Define fundamental weights $\Lambda_i\in (\he)^*$ by $\langle
\Lambda_i,\alpha_j \rangle=\delta_{ij}$ and $\Lambda_i(d)=0$,
$i,j=0,\dots,\ell$.

Let $V$ be a highest weight module for affine Lie algebra $\gt$.
Then $V$ is generated by a highest weight vector $v_\Lambda$ such
that
\begin{eqnarray*}
h\cdot v_\Lambda & = & \Lambda(h) v_\Lambda,\quad\textrm{for}\ h\in\he, \\
x\cdot v_\Lambda & = & 0,\quad\textrm{for}\ x\in\nt_+,
\end{eqnarray*}
for $\Lambda\in (\he)^*$. Module $V$ is a direct sum of weight
subspaces $V_\mu=\{v\in V \,|\,h\cdot V =\mu(h) v \textrm{ za }
h\in\he\},\,\mu\in\he$.

Standard (i.e. integrable highest weight) $\gt$-module $L(\Lambda)$
is an irreducible highest weight module, with the highest weight
$\Lambda$ being dominant integral, i.e.
$$\Lambda=k_0 \Lambda_0+k_1 \Lambda_1+\dots+k_\ell \Lambda_\ell,$$
where $k_i\in\Z_+$, $i=0,\dots,\ell$. The central element $c$ acts on
$L(\Lambda)$ as multiplcation by scalar
$$k=\Lambda(c)=k_0+k_1+\dots+k_\ell,$$
which is called the level of module $L(\Lambda)$.

\section{Feigin-Stoyanovsky's type subspace}

Vector $v\in\h$ is said to be \emph{cominuscule} if
$$\{\alpha(v)\,|\,\alpha \in R\}=\{-1,0,1\}.$$
Similarly, weight $\omega\in P$ is said to be \emph{minuscule} if
$$\{ \langle\omega,\alpha\rangle \,|\,\alpha \in R\}=\{-1,0,1\}.$$
One immediately sees that a dominant integral weight $\omega\in P^+$
is minuscule if and only if
$$\langle\omega,\theta\rangle = 1.$$
So, there exist a finite number of minuscule weights. Furthermore,
a vector $v\in\h$ is cominuscule if and only if it is dual to some
minuscule fundamental weight $\omega$, in the sense that
$$v=\nu^{-1}(\omega),$$
for some choice of positive roots.

Fix a cominuscule vector $v\in\h$. In the case of $\g=\gsl(\ell+1,\C)$, all
 fundamental weights are minuscule. Then we can assume that the
cominuscule vector $v$ is dual to a fundamental weight
$$\omega=\omega_m,$$ for some $m \in \{1,\dots, \ell\}$. 
Set
$$\Gamma =
\{\,\alpha \in R \mid \alpha(v) = 1\}=\{\,\alpha \in R \mid \langle\omega, \alpha\rangle = 1\}.$$ Then we have the induced
$\mathbb Z$-gradation of ${\mathfrak g}$:
\begin{equation}
\label{GDecomp_jed} \mathfrak g  =
\mathfrak g_{-1} \oplus \mathfrak g_0 \oplus \mathfrak g_1, \end{equation}
where
\begin{eqnarray*} {\mathfrak g}_0 & = & {\mathfrak h} \oplus
\sum_{\alpha(v)=0}\, {\mathfrak g}_\alpha
\\
 \displaystyle {\mathfrak g}_{\pm1} & = &
\sum_{\alpha \in \pm \Gamma}\, {\mathfrak g}_\alpha.
\end{eqnarray*}
Subalgebras ${\mathfrak g}_1$ and ${\mathfrak g}_{-1}$ are
commutative, and ${\mathfrak g}_0$ acts on them by adjoint action.
The subalgebra ${\mathfrak g}_0$ is reductive with semisimple part
${\mathfrak l}_0=[{\mathfrak g}_0,{\mathfrak g}_0]$ of the type
$A_{m-1}\times A_{\ell-m}$; as a root basis one can take
$\{\alpha_1,\dots,\alpha_{m-1}\}\cup\{\alpha_{m+1},\dots,\alpha_\ell\}$,
and the center is equal to $\C v$.

We illustrate decomposition \eqref{GDecomp_jed} on the picture \ref{ZGrad_fig}, which corresponds to the usual realization of $\g$ as matrices of trace $0$. In this case the subalgebra ${\mathfrak g}_0$ consists of block-diagonal matrices, while ${\mathfrak g}_1$ and $ {\mathfrak g}_{-1}$ consist of
 matrices with non-zero entries only in the upper right or lower-left block, respectively.
\begin{figure}[ht] \caption{}\label{ZGrad_fig}
\begin{center}\begin{picture}(145,140)(-15,-10)
\thicklines \put(0,0){\line(0,1){120}} \put(0,0){\line(1,0){120}}
\put(48,0){\line(0,1){120}} \put(0,72){\line(1,0){120}}
\put(120,120){\line(0,-1){120}} \put(120,120){\line(-1,0){120}}
\put(-15,110){$\g$:} \put(21,93){$\g_0$} \put(81,33){$\g_0$}
\put(18,33){$\g_{-1}$} \put(81,93){$\g_1$}  \linethickness{.075mm}
\multiput(0,76)(0,4){4}{\line(1,0){48}}
\multiput(48,4)(0,4){7}{\line(1,0){72}}
\multiput(52,72)(4,0){8}{\line(0,1){48}}
\multiput(0,100)(0,4){5}{\line(1,0){48}}
\multiput(48,40)(0,4){9}{\line(1,0){72}}
\multiput(92,72)(4,0){8}{\line(0,1){48}}
\multiput(0,92)(0,4){2}{\line(1,0){19}}
\multiput(48,92)(0,4){2}{\line(-1,0){15}}
\multiput(48,32)(0,4){2}{\line(1,0){31}}
\multiput(120,32)(0,4){2}{\line(-1,0){27}}
\multiput(84,72)(4,0){2}{\line(0,1){16}}
\multiput(84,120)(4,0){2}{\line(0,-1){18}}
\end{picture}\end{center}
\end{figure}

Basis of a subalgebra ${\mathfrak g}_1$ can be identified with the set
of roots $\Gamma$. We will call elements  $\gamma\in\Gamma$ {\em
colors} and the set $\Gamma$ {\em the  set of colors}.  In the case
of $\g=\gsl(\ell+1,\C)$, $\omega=\omega_m$, the set of colors is
$$\Gamma=\{\gamma_{ij}\,|\, i=1,\ldots,m; j=m,\ldots,\ell\}$$
where
\begin{equation}\label{GammaIJ_jed} \gamma_{ij}=\alpha_i+\cdots+\alpha_m+\cdots+\alpha_j.\end{equation}
The maximal root $\theta$ is equal to $\gamma_{1\ell}$.

We picture the
set of colors $\Gamma$ as a rectangle with
row-indices $1,\dots,m$ and column-indices $m,\dots,\ell$ (see figure \ref{Gamma_fig}).

Similarly, one also has the induced $\Z$-gradation of  affine Lie
algebra $\gt$:
\begin{eqnarray*}
\gt_0 & = & {\mathfrak g}_0\otimes\C [t,t^{-1}]\oplus \C c \oplus \C d,\\
\gt_{\pm 1} & = & {\mathfrak g}_{\pm 1}\otimes\C [t,t^{-1}],\\
\gt & = & \gt_{-1} + \gt_0 + \gt_1.
\end{eqnarray*}
As above, $\gt_{-1}$ and $\gt_1$ are commutative subalgebras, and
$\gt_1$ is a $\gt_0$-module.

For a dominant integral weight $\Lambda$, we define a
\emph{Feigin-Stoyanovsky's type subspace}
$$W(\Lambda)=U(\gt_1)\cdot v_\Lambda\subset L(\Lambda).$$

Our objective is to find a combinatorial basis of  $W(\Lambda)$. Set
$$\gt_1^+=\gt_1\cap \nt_+,\, \gt_1^-=\gt_1\cap \nt_-.$$ Then we have
$$W(\Lambda)=U(\gt_1^-)\cdot v_\Lambda.$$
By Poincar\'e-Birkhoff-Witt theorem, we have a spanning set of
$W(\Lambda)$ consisting of monomial vectors
\begin{equation}\label{PBWgen_jed}\{x_{\gamma_1}(-n_1)
x_{\gamma_2}(-n_2)\cdots x_{\gamma_r}(-n_r) v_\Lambda\,|\, r\in\Z_+;
\gamma_j\in \Gamma, n_j\in\N\}.\end{equation}

In the end, we'll say a few words about notation.  Elements of the
spanning set \ref{PBWgen_jed} can be identified with monomials from
$U(\gt_1)=S(\gt_1)$. With this in mind, we'll often refer to
elements of $\{x_\gamma(-j) \mid \gamma\in\Gamma,j\in \Z\}$ in
$\gt_1$ as to {\em variables}, {\em elements} or {\em factors}
of a monomial.

We can also identify monomials from $S(\gt_1)$ with {\em colored
partitions}. From the beginnings of the representation-theoretic
approach to Rogers-Ramanujan identities, combinatorial basis of
certain representations were parameterized by partitions
satisfying certain conditions (cf. [LW], [LP]). Let
$\pi:\{x_\gamma(-j) \mid \gamma\in\Gamma,j\in \Z\}\to \Z_+$ be a
colored partition (cf. [P1], section 3). The corresponding monomial
$x(\pi)\in S(\gt_1)$ is
$$x(\pi)=x_{\gamma_1}(-j_1)^{\pi(x_{\gamma_1}(-j_1))}\cdots x_{\gamma_t}(-j_t)^{\pi(x_{\gamma_t}(-j_t))}.$$
From this identification we'll take notation $x(\pi)$ for the
monomials from $S(\gt_1)$. It will be convenient to define some new monomials
by using this identification. Also, our combinatorial conditions for the basis elements can be written in terms of exponents $\pi(x_\gamma(-j))$, 
which gives a parametrization of the basis by a certain generalization of $(k,\ell +1)$-admissible configurations from [FJLMM].
This will prove to be useful in a higher-level case (cf. [T]).

\section{Order on the set of monomials}

\label{uredjaj_sect}

We introduce a linear order on the set of monomials.

On the weight and root lattice, we have an order $\prec$ defined
in the standard fashion: for $\mu,\nu\in P$ set $\mu\prec\nu$ if $\mu-\nu$
is an integral linear combination of simple roots $\alpha_i, i=1,\dots,\ell$, with non-negative coefficients.

Next, we define a linear order $<$ on the set of colors $\Gamma$ which is an extension of the order $\prec$.
For elements of $\Gamma$, $\gamma_{i'j'}\prec\gamma_{ij}$ is equivalent to saying that
$i'\geq i$ and $j'\leq j$. The order $<$ on $\Gamma$ is defined in the following way:
$$\gamma_{i'j'}<\gamma_{ij}\textrm{\quad if \quad}\left\{\begin{array}{l}i'>i \\ i'=i,\ j'<j.  \end{array}\right.$$
It is clear that this is a linear order on the set of colors.

On the set of variables
$\{x_\gamma(-n)\,|\,\gamma\in\Gamma,\,n\in\Z\}\subset \gt_1$ we
define a linear order $<$ so that we compare degrees first, and then colors of
variables:
$$x_\alpha(-i)<x_\beta(-j) \textrm{\quad if \quad}
\left\{\begin{array}{l}
-i<-j,\\
i=j \textrm{\quad and \quad} \alpha<\beta.
\end{array}\right. $$

Since the algebra $\gt_1$ is commutative, we can assume that the
variables in monomials from $S(\gt_1)$  are sorted ascendingly from
left to right. The order $<$ on the set of monomials is defined as a
lexicographic order, where we compare variables from right to left
(from the greatest to the lowest one). If $x(\pi)$ and $x(\pi')$ are
two monomials,
\begin{eqnarray*}
x(\pi) & = & x_{\gamma_r}(-n_r) x_{\gamma_{r-1}}(-n_{r-1})\cdots
x_{\gamma_2}(-n_2)x_{\gamma_1}(-n_1),\\
x(\pi') & = & x_{\gamma_s'}(-n_s')
x_{\gamma_{s-1}'}(-n_{s-1}')\cdots
x_{\gamma_2'}(-n_2')x_{\gamma_1'}(-n_1'),
\end{eqnarray*}
then $x(\pi)<x(\pi')$ if there exist $i_0\in \N$ so that
$x_{\gamma_i}(-n_i)=x_{\gamma_i'}(-n_i'),\ \textrm{for all}\ i <
i_0,$ and either $i_0=r+1\leq s$ or
$x_{\gamma_{i_0}}(-n_{i_0})<x_{\gamma_{i_0}'}(-n_{i_0}')$.

This monomial order is compatible with multiplication:
\begin{prop} \label{uredjaj}
Let $$x(\pi_1)\leq x(\mu_1)\quad \textrm{and} \quad x(\pi_2) \leq
x(\mu_2).$$ Then $$x(\pi_1)x(\pi_2) \leq x(\mu_1)x(\mu_2),$$ and if
one of the first two inequalities is strict, then the last one is
also strict.
\end{prop}

\begin{dokaz}
 By the definition of the order $<$, we compare two monomials so that we compare their greatest elements first. Let
$x_{\alpha_1}(-j_1)$, $x_{\alpha_2}(-j_2)$, $x_{\beta_1}(-i_1)$,
$x_{\beta_2}(-i_2)$ be the greatest variables in $x(\pi_1)$,
$x(\pi_2)$, $x(\mu_1)$, $x(\mu_2)$ respectively. Then
$x_{\alpha_1}(-j_1)\leq x_{\beta_1}(-i_1)$ and
$x_{\alpha_2}(-j_2)\leq x_{\beta_2}(-i_2)$. The greatest element in
$x(\pi)$ we'll be greater of the two $x_{\alpha_1}(-j_1)$ and
$x_{\alpha_2}(-j_2)$; one can assume it to be $x_{\alpha_1}(-j_1)$.
Similarly, greatest element in $x(\mu)$ we'll be greater of the two
 $x_{\beta_1}(-i_1)$ and $x_{\beta_2}(-i_2)$. There are two
possibilities:
\begin{enumerate}
\renewcommand{\labelenumi}
{(\roman{enumi})}
\item the greatest element of $x(\mu)$ is strictly greater than the greatest element of $x(\pi)$. In that case
$x(\pi)<x(\mu)$.
\item the greatest element of $x(\mu)$ is equal to the greatest element of $x(\pi)$.
Then $x_{\alpha_1}(-j_1)=x_{\beta_1}(-i_1)$ and we can take
$x_{\beta_1}(-i_1)$ for the greatest element of $x(\mu)$. We proceed
by induction: let $x(\pi_1')$ and $x(\mu_1')$ be monomials gotten
from $x(\pi_1)$ and $x(\mu_1)$, respectively, by omitting
$x_{\alpha_1}(-j_1)=x_{\beta_1}(-i_1)$. Then $x(\pi_1')\leq
x(\mu_1')$ and we can continue to apply the same procedure to
monomials $x(\pi_1')$, $x(\pi_2)$, $x(\mu_1')$ and $x(\mu_2)$. After
a finite number of steps either case $(i)$ will occur, or we'll
exhaust monomials $x(\pi_1)$ and $x(\pi_2)$. Both these cases imply
$x(\pi)\leq x(\mu)$, and the equality occurs only if both initial
inequalities were in fact equalities.
\end{enumerate}
\end{dokaz}

For monomials from $S(\gt_1)$, we also define degree
and shape of a monomial. \emph{Degree} of a monomial is equal to the
sum of degrees of its variables. For
$$x(\pi) = x_{\gamma_r}(-n_r) x_{\gamma_{r-1}}(-n_{r-1})\cdots
x_{\gamma_2}(-n_2)x_{\gamma_1}(-n_1),$$ its degree is equal to
$-n_1-n_2-\dots-n_r$. A \emph{shape} of a monomial is gotten from its colored
partition by forgetting colors and considering only degrees of
factors. More precisely, for a monomial $x(\pi)$ and its partition
$\pi:\{x_\gamma(-n) \mid \gamma\in\Gamma,n\in \Z\}\to \Z_+$, the
corresponding shape will be
\begin{eqnarray*}
& & s_\pi:\Z\to\Z_+,\\
& & s_\pi(j)=\sum_{\gamma\in\Gamma}\pi(x_\gamma(-j)).
\end{eqnarray*}
A linear order can also be defined on the set of shapes; we'll say
that $s_\pi<s_{\pi'}$ if there exists $j_0\in\Z$ such that
$s_\pi(j)=s_{\pi'}(j)$ for  $j<j_0$ and either $s_\pi(j_0)<
s_{\pi'}(j_0)$ and $s_\pi(j')\neq 0$ for some $j'>j_0$, or
$s_\pi(j_0)>s_{\pi'}(j_0)$ and $s_\pi(j)=0$ for $j>j_0$.

In the end, for the sake of simplicity, we introduce the following
notation:
$$x_{rs}(-j)=x_{\gamma_{rs}}(-j),$$
 for $\gamma_{rs}\in\Gamma$.

\section{Vertex operator construction of level 1 modules}

\label{voakonstr_sect} We use the vertex operator algebra construction of the basic $\gt$-modules (i.e. the standard
$\gt$-modules of level
$1$). We'll sketch this construction in this section, details can be found in [FLM], [DL] or [LL]; see also [FK],
[S].

We have denoted by $P$ and $Q$ weight and root lattices of
$\g$, respectively. There exists a central extension $\hat{P}$ of $P$ by the finite cyclic group $\langle e^{\pi i/(\ell+1)^2}\rangle $ of
order $2(\ell+1)^2$,
$$1\longrightarrow \langle e^{\pi i/(\ell+1)^2}\rangle \longrightarrow \hat{P}\longrightarrow P \longrightarrow 1.$$
By restricting, one gets a central extension
$\hat{Q}$ of $Q$. Central extension can be chosen such that the
corresponding 2-cocycle 
$$\epsilon : P\times P \to \langle e^{\pi i/(\ell+1)^2}\rangle $$
satisfies
$$\epsilon(\alpha,\beta)/\epsilon(\beta,\alpha)=(-1)^{\langle \alpha,\beta\rangle}\quad \textrm{for}\ \alpha,\beta\in Q.$$
Let
$$c(\lambda,\mu)=\epsilon(\lambda,\mu)/\epsilon(\mu,\lambda)\quad \textrm{for}\
\lambda, \mu \in P$$ be the corresponding bimultiplicative,
alternating commutator map (cf. [FLM]).

Inside $\gt$ there is a Heisenberg subalgebra
$$\hz=\sum_{n\in\Z \setminus\{0\}}\h\otimes t^n \oplus \C c.$$
 We also introduce subalgebras
\begin{eqnarray*}
\hk & = & \h\otimes \C [t,t^{-1}]\oplus \C c,\\
\hk_\pm & = & \h\otimes t^{\pm 1}\C [t^{\pm 1}],
\end{eqnarray*}
and by $\C [P]$ and $\C [Q]$ we denote group algebras of weight and
root lattices, respectively. Bases of $\C [P]$ and $\C [Q]$ consist of
elements $\{e^\lambda\,|\,\lambda\in P\}$ and
$\{e^\alpha\,|\,\alpha\in Q\}$, respectively.

Consider the induced $\hz$-module
$$M(1)=U(\hk)\otimes_{\h\otimes \C[t]\oplus \C c}\C,$$
where $\h\otimes \C[t]$ acts trivially on $\C$, and $c$ acts as $1$.
Module $M(1)$ is irreducible module for the Heisenberg subalgebra
$\hz$; as a vector space, $M(1)$ is naturally isomorphic to the
symmetric algebra $S(\hk_-)$ (cf. [FLM]).

Consider tensor products
\begin{eqnarray*}
V_P & = & M(1)\otimes \C [P],\\
V_Q & = & M(1)\otimes \C [Q];
\end{eqnarray*}
there is a natural inclusion $V_Q\subset V_P$. For simplicity, we
will often write $e^\lambda$ instead of $1\otimes e^\lambda$, and
$1$ instead of $1\otimes 1$.

Space $V_P$ carries a $\hk$-module structure: $\hz$ acts as
$\hz\otimes 1$ and $\h\otimes t^0$ acts as $1\otimes \h$. Operators
$h(0),h\in \h$ on $\C [P]$ are defined as follows
$$h(0)\cdot e^\lambda=\langle h,\lambda \rangle e^\lambda$$
for $\lambda\in P$. On $V_P$ we have also the action of
the group algebra $\C [P]$:
$$e^\lambda = 1\otimes e^\lambda,\quad \lambda \in P,$$
where the latter operator $e^\lambda$ is a multiplication in $\C
[P]$. It will be clear from the context when $e^\lambda$ represents
a multiplication operator, and when an element of $V_P$. Define also
operators $\epsilon_\lambda$ by
$$\epsilon_\lambda \cdot e^\mu= \epsilon(\lambda,\mu) e^\mu,$$
 for $\lambda,\mu \in P$.

For elements of $V_P$ define a degree: for $v=h_1(-n_1)\cdots
h_r(-n_r)\otimes e^\lambda$ set
$$\deg(v)=-n_1-n_2-\dots -n_r - \frac{1}{2}\langle \lambda,\lambda
\rangle.$$ This gives a grading on $V_P$, which is bounded from
above.

We will use independent commuting formal variables
$z,z_0,z_1,z_2,\dots$. For a vector space $V$, denote by $V [[z]]$
the space of formal series of nonnegative integral powers of $z$
with coefficients in $V$. Similarly, denote by $V[[z,z^{-1}]]$ the
space of formal Laurent series, and by $V\{z\}$ the space of formal
series of rational powers of $z$ with coefficients in $V$.

Define also one more family of operators $z^h\in
 (\textrm{End}\,V_P)\{z\}$ by
$$z^h \cdot e^\lambda=e^\lambda z^{\langle h,\lambda \rangle},$$
for $h\in \h, \lambda\in P$.

Space $V_Q$ has a natural structure of vertex operator algebra and
$V_P$ is a module for this algebra (cf. [FLM],[DL]). Before we define VOA-structure
on $V_Q$, define operators
\begin{eqnarray*}
h(z) & = & \sum_{j\in\Z}h(j) z^{-j-1},\\
E^{\pm}(h,z)& =& \exp \left(\sum_{m\geq 1}h(\pm m) \frac{z^{\mp m}}{\pm m}\right),
\end{eqnarray*}
for $h\in \h$. We define vertex operators for all elements of $V_P$,
rather than just for elements of $V_Q$. For the lattice elements,
i.e. for the elements $1\otimes e^\lambda=e^\lambda$ set:
\begin{equation}Y(e^\lambda,z)=E^-(-\lambda,z)E^+(-\lambda,z)\otimes
e^\lambda z^\lambda \epsilon_\lambda.\label{vop_jed}\end{equation}
Generally, for a homogenous vector $v\in V_P$
$$v=h_1(-n_1)\cdots h_r(-n_r)\otimes e^\lambda,$$
 $n_1,\dots,n_r\geq 1$, set
$$Y(v,z)=\nop\left(\frac{\partial_z^{n_1-1}}{(n_1-1)!}h_1(z)\right)\cdots \left(\frac{\partial_z^{n_r-1}}{(n_r-1)!}h_r(z)\right) Y(
e^\lambda, z)\nop,$$ where  $\nop\cdot\nop$ is a normal ordering
procedure (cf. [FLM]), meaning that coefficients in the enclosed expression should be rearranged
in a way that in every product all the operators $h(m), h\in\h,m<0$ are placed to the
left of the operators $h(m),h\in\h,m\geq 0$. This way we get a well
defined linear map
\begin{eqnarray*}
Y:V_P & \to & (\End V_P)\{z\},\\
v & \mapsto & Y(v,z).
\end{eqnarray*}

By using vertex operators, one can define a structure of
$\gt$-module on $V_P$. For $\alpha\in R$ set
$$x_\alpha(z)=\sum_{j\in\Z}x_\alpha (j) z^{-j-1}=Y(e^\alpha,z).$$
Actions of $h(j)$ and $c$ have already been defined, and $d$ acts as a
degree operator. Then the cosets $V_Q$ and $V_Q
e^{\omega_j},j=1,\dots,\ell$ become standard $\gt$-modules of level
$1$ with highest weight vectors $v_0=1$ and
$v_j=e^{\omega_j},j=1,\dots,\ell$, respectively (cf. [FLM,DL]).  Moreover,
$$L(\Lambda_0)\cong V_Q,\ L(\Lambda_j)\cong V_Q
e^{\omega_j} \textrm{ for }j=1,\dots,\ell$$ and
$$V_P\cong L(\Lambda_0)\oplus L(\Lambda_1)\oplus\dots
L(\Lambda_\ell).$$


Vertex operators $Y(v,z),v\in V_P$ satisfy (generalized) Jacobi
identity. It will be of importance to us a variant of that identity
in the case when vectors $u,v$ are of type $u=u^*\otimes
e^\lambda,v=v^*\otimes e^\mu$, for $\lambda\in Q,\mu \in P,\
u^*,v^*\in M(1)$, or, even more special, when $u=1\otimes
e^\lambda,v=1\otimes e^\mu$, for $\lambda\in Q,\mu \in P$. Then one
has
\begin{eqnarray*}  & \hspace{-2ex} z_{0}^{-1} \delta \left (\frac{z_{1} - z_{2}}{z_{0}} \right
)\hspace{-.5ex} Y(u,z_{1})Y(v,z_{2}) \hspace{-.3ex}-\hspace{-.4ex}
(-1)^{\langle\lambda, \mu\rangle}c(\lambda ,\mu) z_{0}^{-1} \delta
\left (\frac{z_{2} - z_{1}}{-z_{0}} \right
)\hspace{-.5ex}Y(v,z_{2})Y(u,z_{1}) = & \nonumber\\& =\; z_{2}^{-1}
\delta \left ( \frac{z_{1} - z_{0}}{z_{2}} \right ) Y(Y(u,z_{0})v,
z_{2}),&\end{eqnarray*} where $\delta(z)=\sum_{n\in \Z} z^n$ is a
formal delta-function (cf. [FLM],[LL]), and binomial expressions that appear in
expansions of delta-functions are understood to be expanded in
nonnegative terms of the second variable.

Next we introduce intertwining operators $\mathcal Y$. For $\mu\in
P,\,v=v^*\otimes e^\mu$ define $${\mathcal Y}
(v,z)=Y(v,z)e^{i\pi\mu}c(\cdot,\mu).$$ This way we obtain a map
\begin{eqnarray*}
{\mathcal Y}:V_P & \to & (\End V_P)\{z\},\\
 v & \mapsto & {\mathcal Y}
(v,z).\end{eqnarray*} Then we have (ordinary) Jacobi identity
\begin{eqnarray*}
&z_{0}^{-1} \delta \left (\frac{z_{1} - z_{2}}{z_{0}} \right
)Y(u,z_{1}){\mathcal Y}(v,z_{2}) \;- z_{0}^{-1} \delta \left
(\frac{z_{2} - z_{1}}{-z_{0}} \right ){\mathcal Y}(v,z_{2})Y(u,z_{1}) =
\\& = z_{2}^{-1} \delta \left ( \frac{z_{1} - z_{0}}{z_{2}} \right )
{\mathcal Y}(Y(u,z_{0}) v , z_{2}).\end{eqnarray*} For  $\mu\in Q$,
operators ${\mathcal Y} (v,z)$ are equal to vertex operators
$Y(v,z)$. Restrictions of ${\mathcal Y} (v,z)$ are in fact maps
\begin{equation} \label{evop_jed} {\mathcal Y} (v,z):L(\Lambda_i)\to
L(\Lambda_j)\{z\},\end{equation} if $\mu+\omega_i \equiv \omega_j
\mod Q$. So, restrictions of $\mathcal Y$ define intertwining
operators between standard modules of level $1$ ([DL]).


Consider now a special case when $v=e^\mu$. It is interesting to
know when the operators ${\mathcal Y}(e^\mu,z_2)$ from
\eqref{evop_jed} commute with the action of $\gt_1$, i.e. when
$$[Y(e^\gamma,z_1),{\mathcal
Y}(e^\mu,z_2)]=0,\quad \gamma\in\Gamma.$$ By the commutator formula
for intertwining operators ([DL]) that is equivalent to
$$Y(e^\gamma,z_0)e^\mu \in V_P[[z_0]],$$
for all $\gamma\in\Gamma$. From the definition of vertex operators
\eqref{vop_jed} one gets
\begin{equation}\label{djel_vop_jed}
\qquad Y(e^\gamma,z_0)e^\nu=C e^{\gamma+\nu} z_0^{\langle \gamma,\nu
\rangle}+\underbrace{\qquad\dots\qquad}_{\textrm{higher power terms}} \quad \in \ \ z_0^{\langle \gamma,\nu
\rangle}V_P[[z_0]],
\end{equation} for some $C\in \C^\times$.
So, operators ${\mathcal Y}(e^\mu,z_2)$  commute with $\gt_1$ if and
only if
$$\langle \gamma,\mu \rangle \geq 0,\quad \textrm{for all}\
\gamma\in\Gamma.$$ In the section \ref{IntOp_sect}, we describe all $\mu\in P$ that satisfy
this relation.

\section{Operator $e(\omega)$}

For $\lambda\in P$, $e^\lambda$ denotes multiplication operator
$1\otimes e^\lambda$ in $V_P=M(1)\otimes \C [P]$. Set
$$e(\lambda)=e^\lambda \epsilon(\cdot,\lambda),\qquad e(\lambda):V_P\to V_P.$$
Clearly, $e(\lambda)$ is a linear bijection. Its restrictions on
standard modules are bijections from one fundamental module
$L(\Lambda_i)$ onto another fundamental module $L(\Lambda_{i'})$. From the
definition of vertex operators \eqref{vop_jed}, one gets
the following commutation relation
$$Y(e^\alpha,z)e(\lambda)=e(\lambda)
z^{\langle\lambda,\alpha\rangle} Y(e^\alpha,z),$$
for $\alpha\in R$. In terms of
components, we have
\begin{equation} x_\alpha(n)
e(\lambda)=e(\lambda) x_\alpha(n+\langle\lambda,\alpha\rangle),\quad
n\in\Z. \label{komut_ea_xn_jed}\end{equation}

For $\lambda=\omega$ and $\gamma\in \Gamma$, the relation
\eqref{komut_ea_xn_jed} becomes
$$x_\gamma(n) e(\omega)=e(\omega)
x_\gamma(n+1).$$

More generally, for a monomial $x(\pi)\in S(\gt_1)$, denote by
$x(\pi^+)\in S(\gt_1)$ the monomial corresponding to the partition $\pi^+$,
defined by $\pi^+(x_\gamma(n+1))=\pi(x_\gamma(n))$. We can say that $x(\pi^+)$
is obtained from $x(\pi)$ by raising degrees of all its factors by
$1$. Then $$x(\pi)e(\omega)=e(\omega) x(\pi^+).$$

\section{Difference and initial conditions}

Initial conditions for the level $1$ standard module $L(\Lambda_i)$ are consequence
of a simple observation that monomials from the monomial basis can't
contain elements of degree $-1$ that act as zero on the highest weight
vector $v_i$ of $L(\Lambda_i)$. So,
 we have to establish for which $\gamma\in\Gamma$,
elements $x_\gamma(-1)$ annihilate $v_i$. Then we can exclude from
the spanning set \eqref{PBWgen_jed} all monomials $x(\pi)$ that
contain such factors.

Since $v_i=e^{\omega_i}$, for $i=1,\dots,\ell$, and $v_0=1=e^0$, relation \eqref{djel_vop_jed} gives
\begin{equation} \label{GDjel0_jed}
x_\gamma(z)v_i=\left(\sum_{j\in\Z}x_\gamma(-j) z^{j-1} \right)v_i \in z^{\langle \gamma,\omega_i \rangle} (v_i V_Q)[[z]],\end{equation}
and
\begin{equation} \label{GDjel_jed}
x_\gamma(z)v_i=\left(\sum_{j\in\Z}x_\gamma(-j) z^{j-1} \right)v_i \in z^{\langle \gamma,\omega_i \rangle} (v_i V_Q)[[z]],\end{equation}
for $i=1,\dots,\ell$.
Since by \eqref{GammaIJ_jed}, $\langle \gamma_{rs},\omega_i\rangle=1$ if $r\leq i\leq s$, and zero otherwise, by comparing constant terms in \eqref{GDjel0_jed} and \eqref{GDjel_jed}, we get
\begin{equation} \label{GDjelIC_jed}
x_{\gamma_{rs}}(-1)v_i=\left\{\begin{array}{l l}
0, & r\leq i\leq s, \\
C e^{\gamma_{rs}},\,C\in\C^\times, & i=0,\\
C e^{\gamma_{rs}+\omega_i},\,C\in\C^\times, & \textrm{otherwise}.
\end{array}\right.\end{equation}

For a monomial $x(\pi)\in S(\gt_1^-)$ we say that it satisfies
\emph{initial conditions} for $W(\Lambda_i)$, if it doesn't contain
factors of degree $-1$ that annihilate $v_i$. We'll often abbreviate
this by saying that $x(\pi)$ satisfies IC for $W(\Lambda_i)$. From \eqref{GDjelIC_jed} we see that $x(\pi)$ satisfies initial conditions on $W(\Lambda_i)$
if the colors of elements of degree $-1$ lie below the  $i$-th row
(in case $i\leq m$), or, to the left of the $i$-th column (for
$i\geq m$). 

%
%
%
%
%
%
%

Difference conditions will be consequences of relations between
operators $x_\gamma(z)$, and fortiori, between monomial vectors
$x(\pi)v_i$.

To obtain these, consider the basic module
$L(\Lambda_0)$ with highest weight vector $v_0=1=e^0$ (cf. section \ref{voakonstr_sect}).
This is a vertex operator algebra, with $1$ as the vacuum element, and $L(\Lambda_i)$ is
a module for this algebra.
We are looking for relations between vectors of type $$x_\gamma(-1) x_{\gamma'}(-1)1,\quad \gamma,\gamma'\in\Gamma.$$
These will in turn induce relations between corresponding vertex operators on $L(\Lambda_i)$.

From \eqref{GDjel0_jed} we have
$$x_\gamma(-1) x_{\gamma'}(-1)1=x_\gamma(-1) e^{\gamma'}.$$
Since
$$\langle \gamma,\gamma' \rangle= \left\{\begin{array}{l l}
2, & \gamma=\gamma',\\
1, & \textrm{$\gamma$ and $\gamma'$ lie in the same row or column}, \\
0 & \textrm{otherwise},
\end{array}\right.$$
relation \eqref{djel_vop_jed} implies
$$x_\gamma(-1) x_{\gamma'}(-1)1=\left\{\begin{array}{l l}
0, & \textrm{\begin{minipage}[t]{4cm}$\gamma$ and $\gamma'$ lie in the same row or column,\end{minipage}} \\
Ce^{\gamma+\gamma'},\,C\in\C^\times, & \textrm{otherwise}.
\end{array}\right.$$
Fix two rows $r_1<r_2$ and two columns $s_2<s_1$. Note that
$$\gamma_{r_1 s_1}+\gamma_{r_2 s_2}=\gamma_{r_1 s_2}+\gamma_{r_2
s_1}.$$ This gives us
$$x_{r_2 s_2}(-1) x_{r_1 s_1}(-1)1=C\cdot x_{r_2 s_1}(-1) x_{r_1
s_2}(-1)1,$$ for some $C\in\C^\times$.

We've obtained two types of relations:
$$x_\gamma(-1) x_{\gamma'}(-1)1=0,$$
if $\gamma$ and $\gamma'$ lie in the same row/column, and
$$x_\gamma(-1) x_{\gamma'}(-1)1=C\cdot x_{\gamma_1}(-1)
x_{\gamma_1'}(-1)1,$$ if $\gamma$, $\gamma_1$, $\gamma'$  i
$\gamma_1'$ are vertices of a rectangle in $\Gamma$, as in the figure \ref{Rel_fig}.

\begin{figure}[ht] \caption{} \label{Rel_fig}
\begin{center}\begin{picture}(200,140)(-8,-10) \thicklines
\put(0,0){\line(1,0){180}} \put(0,0){\line(0,1){120}}
\put(180,120){\line(0,-1){120}} \put(180,120){\line(-1,0){180}}
\put(-6,111){$\scriptstyle 1$} \put(-6,99){$\scriptstyle 2$}
\put(-8,3){$\scriptstyle m$} \put(3,-8){$\scriptstyle m$}
\put(174,-8){$\scriptstyle \ell$}

\linethickness{.075mm} \multiput(0,40)(4,0){35}{\line(1,0){2}}
\multiput(0,80)(4,0){35}{\line(1,0){2}}
\multiput(90,0)(0,4){20}{\line(0,1){2}}
\multiput(140,0)(0,4){20}{\line(0,1){2}}

\put(-8,37){$\scriptstyle r_2$}  \put(-8,77){$\scriptstyle r_1$}
\put(87,-8){$\scriptstyle s_2$} \put(137,-8){$\scriptstyle s_1$}

\put(88,38){$\scriptscriptstyle \bullet$}
\put(81,36){$\scriptstyle \gamma$}
\put(139,38){$\scriptscriptstyle \bullet$}
\put(143,36){$\scriptstyle \gamma_1$}
\put(88,78){$\scriptscriptstyle \bullet$}
\put(81,80){$\scriptstyle \gamma_1'$}
\put(139,78){$\scriptscriptstyle \bullet$}
\put(143,80){$\scriptstyle \gamma'$}
\end{picture}\end{center}
\end{figure}

Since the algebra $\gt_1$ is commutative, vertex operators
$Y(x_\gamma(-1)x_{\gamma'}(-1)1,z)$ are equal to products of $x_{\gamma}(z)$ and $x_{\gamma'}(z)$ as ordinary products of Laurent series  (cf. [DL],[LL]).
This way we get relations between vertex operators on level $1$ modules
\begin{eqnarray}
\label{rel1_jed} x_\gamma(z) x_{\gamma'}(z) & = & 0,\\
\label{rel2_jed}  x_\gamma(z) x_{\gamma'}(z) & = & C\cdot
x_{\gamma_1}(z) x_{\gamma_1'}(z).
\end{eqnarray}

Fix $n\in \N$ and consider the coefficients of $z^{n-2}$ in
\eqref{rel1_jed} and \eqref{rel2_jed}. From the first relation we
have
$$
0 = \sum_{i+j=n}x_\gamma(-i) x_{\gamma'}(-j).$$ 
In each such sum we
can identify the minimal monomial with regard to the ordering $<$,
which is then called \emph{the leading term} of the relation. This can be expressed
in terms of other monomials in the sum, so we can exclude from the
spanning set \eqref{PBWgen_jed} all monomials that contain leading
terms (cf. [P1],[P2]). All the monomials appearing in the sum are of
length $2$ and of total degree $-n$. Because of this, the minimal
between them has to be of the ``minimal shape'', i.e. its factors have
to be either of the same degree (for $n$ even), or degrees have to
differ only for $1$ (for $n$ odd). In the case of even $n$, there is only
one monomial of minimal shape,
\begin{equation} x_\gamma(-j)x_{\gamma'}(-j),\label{vod_clan0_jed}\end{equation}
and that's the leading term of the sum above. For $n$ odd, there
are two monomials of minimal shape,
$$x_{\gamma'}(-j-1)x_{\gamma}(-j),\,
 x_\gamma(-j-1)x_{\gamma'}(-j).$$  By the definition of the order $<$, next we compare colors of elements. First we
 compare colors of elements of degree $-j$, and then of elements of degree $-j-1$.
 If we assume
$\gamma<\gamma'$, then the leading term will be
\begin{equation}x_{\gamma'}(-j-1)x_{\gamma}(-j).\label{vod_clan1_jed}\end{equation}

Analogously we consider relation \eqref{rel2_jed}; we get $$0 =
\sum_{i+j=n}x_\gamma(-i) x_{\gamma'}(-j)-C x_{\gamma_1}(-i)
x_{\gamma_1'}(-j).$$ Assume $\gamma<\gamma_1<\gamma_1'<\gamma'$, as in the figure \ref{Rel_fig}. For
$n$ even we have two monomials of minimal shape
$$x_\gamma(-j)x_{\gamma'}(-j),\,x_{\gamma_1}(-j)x_{\gamma_1'}(-j),$$
and for $n $ odd we have four of them
\begin{eqnarray*} & x_{\gamma'}(-j-1)x_{\gamma}(-j),\,
x_\gamma(-j-1)x_{\gamma'}(-j), &
\\ & x_{\gamma_1'}(-j-1)x_{\gamma_1}(-j),\,
x_{\gamma_1}(-j-1)x_{\gamma_1'}(-j).& \end{eqnarray*} The leading
terms are
\begin{equation}\label{vod_clan2_jed} x_{\gamma_1}(-j)x_{\gamma_1'}(-j)\end{equation} for $n$
even, and
\begin{equation}\label{vod_clan3_jed} x_{\gamma'}(-j-1)x_{\gamma}(-j)\end{equation}
for $n$ odd.

We say that a monomial $x(\pi)\in S(\gt_1^-)$  satisfies
\emph{difference conditions} if it doesn't contain any of the
leading terms \eqref{vod_clan1_jed}, \eqref{vod_clan2_jed},
\eqref{vod_clan3_jed}.

Then, by using proposition \ref{uredjaj}, we get the following proposition
(cf. [P1,Lemma 9.4] and [P2,Theorem 5.3])
\begin{prop}
The set \begin{equation}\label{baza_jed}\{x(\pi)v_i \,|\, x(\pi)
\textrm{ satisfies IC and DC for } W(\Lambda_i)\}\end{equation}
spans $W(\Lambda_i)$.
\end{prop}

Finally, let's have a closer look at the structure of monomials that
satisfy difference and initial conditions for the standard module
$L(\Lambda_i)$ of level 1. Assume that a monomial $x(\pi)$ contains elements
$x_{rs}(-j)$ and $x_{r's'}(-j)$, and $\gamma_{r's'}\leq
\gamma_{rs}$. Then by \eqref{vod_clan0_jed}, $\gamma_{r's'}$ and
$\gamma_{rs}$ cannot lie in the same column or row, because otherwise
$x(\pi)$ would contain a leading term. Hence $\gamma_{r's'}$ and
$\gamma_{rs}$ are opposite vertices of a rectangle in $\Gamma$. By
\eqref{vod_clan2_jed}, they have to be  upper-right and lower-left
vertices of this rectangle, otherwise $x(\pi)$ would contain a
leading term. Since $\gamma_{r's'}\leq \gamma_{rs}$, we conclude
that $r'>r$ and $s'<s$, i.e. $\gamma_{r' s'}$ must lie in the shaded area as illustrated on the figure \ref{DC1_fig}.

\begin{figure}[ht] \caption{} \label{DC1_fig}
\begin{center}\begin{picture}(200,140)(-8,-10) \thicklines
\put(0,0){\line(1,0){180}} \put(0,0){\line(0,1){120}}
\put(180,120){\line(0,-1){120}} \put(180,120){\line(-1,0){180}}
\put(-6,111){$\scriptstyle 1$} \put(-6,99){$\scriptstyle 2$}
\put(-8,3){$\scriptstyle m$} \put(3,-8){$\scriptstyle m$}
\put(174,-8){$\scriptstyle \ell$}

\linethickness{.075mm} \multiput(0,78)(4,0){15}{\line(1,0){2}}
\multiput(58,0)(0,4){20}{\line(0,1){2}}

\put(60,80){$\scriptscriptstyle \bullet$}
\put(64,82){$\scriptstyle \gamma_{rs}$} \put(-6,80){$\scriptstyle
r$} \put(60,-8){$\scriptstyle  s$}

\multiput(0,3)(0,3){11}{\line(1,0){58}}
\multiput(0,42)(0,3){12}{\line(1,0){58}}
\multiput(0,36)(0,3){2}{\line(1,0){22}}
\multiput(42,36)(0,3){2}{\line(1,0){16}}
 \put(24,36){$\scriptstyle
\gamma_{r's'}$}
\end{picture}\end{center}
\end{figure}

\noindent  Next, assume that a monomial $x(\pi)$ contains elements
$x_{rs}(-j)$ and $x_{r's'}(-j-1)$. Then, by a similar argument as
above, one concludes that $r'>r$ or $s'<s$, which is illustrated on the figure
\ref{DC2_fig}.

\begin{figure}[ht] \caption{} \label{DC2_fig}
\begin{center}\begin{picture}(200,140)(-8,-10) \thicklines
\put(0,0){\line(1,0){180}} \put(0,0){\line(0,1){120}}
\put(180,120){\line(0,-1){120}} \put(180,120){\line(-1,0){180}}
\put(-6,111){$\scriptstyle 1$} \put(-6,99){$\scriptstyle 2$}
\put(-8,3){$\scriptstyle m$} \put(3,-8){$\scriptstyle m$}
\put(174,-8){$\scriptstyle \ell$}

\linethickness{.075mm} \multiput(58,78)(4,0){31}{\line(1,0){2}}
\multiput(58,78)(0,4){11}{\line(0,1){2}}

\put(60,80){$\scriptscriptstyle \bullet$}
\put(64,82){$\scriptstyle \gamma_{rs}$} \put(-6,80){$\scriptstyle
r$} \put(60,-8){$\scriptstyle  s$}

\multiput(0,78)(0,3){14}{\line(1,0){58}}
\multiput(0,3)(0,3){11}{\line(1,0){180}}
\multiput(0,42)(0,3){12}{\line(1,0){180}}
\multiput(0,36)(0,3){2}{\line(1,0){52}}
\multiput(72,36)(0,3){2}{\line(1,0){108}}
 \put(54,36){$\scriptstyle
\gamma_{r's'}$}
\end{picture}\end{center}
\end{figure}

From these observations we conclude that colors of the elements of
the same degree  $-j$ inside $x(\pi)$ make a descending sequence as
pictured on the figure \ref{DCSeq_fig}; appropriate row-indices strictly increase, while
column-indices strictly decrease. Set of colors of elements of
degree  $-j-1$ also form a decreasing sequence which is placed below
or on the left of the minimal color of elements of degree $-j$.

Initial conditions for $W(\Lambda_i)$ imply that the sequence of
colors of elements of degree $-1$ lies below the $i$-th row (if
$0\leq i\leq m$), or on the left of the $i$-th column (for $m \leq i
\leq \ell$) (see figure \ref{ICSeq_fig}).

These considerations also imply the following
\begin{prop} \label{DCTrans_prop}
If $x_\gamma(-j)<x_{\gamma'}(-j')<x_{\gamma''}(-j'')$ are such that
monomials $x_\gamma(-j)x_{\gamma'}(-j')$ and $x_{\gamma'}(-j')x_{\gamma''}(-j'')$ satisfy difference conditions, then so does
$x_{\gamma}(-j)x_{\gamma''}(-j'')$, and consequently $x_\gamma(-j)x_{\gamma'}(-j')x_{\gamma''}(-j'')$.
\end{prop}
Hence, under the assumption that factors in monomials are sorted descendingly from right to left, to see if a monomial satisfies difference conditions, it is enough to check difference conditions on all pairs of successive factors in it.

\section{Intertwining operators}

\label{IntOp_sect}

As we've  already seen in section \ref{voakonstr_sect}, operators
$${\mathcal Y} (e^\lambda,z):L(\Lambda_i)\to
L(\Lambda_{i'})\{z\},$$  commute with the action of $\gt_1$ if and only
if \begin{equation} \label{ekomuvj_jed} \langle \lambda,\gamma
\rangle \geq 0,\end{equation} for all $\gamma\in \Gamma$.

Define ``minimal'' weights that satisfy \eqref{ekomuvj_jed}:
\begin{equation} \label{lambda_def_jed}
\begin{array}{l @{\hspace{2cm}} l}
\lambda_1=\omega_1, & \lambda_m'=\omega_m-\omega_{m+1},\\
\lambda_2=\omega_2-\omega_1, & \lambda_{m+1}'=\omega_{m+1}-\omega_{m+2},\\
\lambda_3=\omega_3-\omega_2, &  \quad \vdots \\
 \quad\vdots & \lambda_{\ell-1}'=\omega_{\ell-1}-\omega_{\ell},\\
\lambda_m=\omega_m-\omega_{m-1}, & \lambda_\ell'=\omega_\ell.
\end{array} \end{equation}
Then relation \eqref{GammaIJ_jed} gives
\begin{equation} \label{ScalProd_jed}
\begin{array}{rcl} \langle\lambda_r, \gamma \rangle & = & \left\{\begin{array}{ll}
1, & \textrm{if}\ \gamma\ \textrm{lies in the}\ r \textrm{-th row},
\\
0, & \textrm{otherwise},\end{array}\right.\\
\langle\lambda_s', \gamma \rangle & = & \left\{\begin{array}{ll} 1,
& \textrm{if}\ \gamma\ \textrm{lies in the }\ s \textrm{-th column},
\\
0, & \textrm{otherwise}.\end{array}\right.\end{array}
\end{equation}
It is obvious that every nonnegative $\Z$-linear combination
$\lambda\in P$ of these weights also satisfies condition
\eqref{ekomuvj_jed} and, consequently, the appropriate intertwining
operator ${\mathcal Y}(e^\lambda,z)$ commutes  with $\gt_1$.  It can
easily be shown that a weight $\lambda\in P$ satisfies \eqref{ekomuvj_jed}
if and only if $\lambda$ can be written in this way. For example
\begin{eqnarray}
\omega_3 & = & \lambda_1+\lambda_2+\lambda_3,\nonumber\\
\omega_r & = & \lambda_r+\lambda_{r-1}+\dots +\lambda_1,\quad
\textrm{for}\ r\leq m,\nonumber\\
\omega_s & = & \lambda_s'+\lambda_{j+1}'+\dots +\lambda_\ell',\quad
\textrm{for}\ s\geq m,\nonumber\\
\omega_m & = & \lambda_m+\lambda_{m-1}+\dots
+\lambda_1=\lambda_m'+\lambda_{m+1}'+\dots +\lambda_\ell'.
\label{omega_lambda_jed}
\end{eqnarray}

In the next section we'll need the following lemma
\begin{lm}
\label{gama_lambda} Let $\gamma_{rs}\in\Gamma$. Then
\begin{equation}\label{gama_omega_jed}\qquad \gamma_{rs}=\lambda_r
+\lambda_s'.\end{equation}
\end{lm}

\begin{dokaz}
By the Cartan matrix of $\g$, we have
\begin{eqnarray*}
\alpha_1 & = & 2\omega_1 -\omega_2, \\
\alpha_j & = & -\omega_{j-1}+2\omega_j -\omega_{j+1};\quad j=2,\dots,\ell-1 \\
\alpha_\ell & = & -\omega_{\ell-1}+2\omega_\ell.
\end{eqnarray*}
The claim now follows from \eqref{GammaIJ_jed} and
\eqref{lambda_def_jed}.
\end{dokaz}

\section{Proof of linear independence}

Write a monomial $x(\pi)\in S(\gt_1^-)$ as a product
$x(\pi)=x(\pi_2) x(\pi_1)$, where $x(\pi_1)$ consists of
elements of degree $-1$, and $x(\pi_2)$ consists of elements of
lower degree. The main technical tool in the proof of linear
independence is the following proposition:

\begin{prop} \label{OpIsp1_prop}
Suppose that a monomial  $x(\pi)$ satisfies difference and initial
conditions for a level $1$ standard module $L(\Lambda_i)$. Then there
exists a coefficient $w(\mu)$ of an intertwining operator ${\mathcal
Y}(e^\mu,z)$
$$w(\mu):L(\Lambda_i) \to L(\Lambda_{i'})$$
for some $i'\in\{0,\dots,\ell\}$, such that:
\begin{itemize}
\item $w(\mu)x(\pi_1)v_i = C e(\omega)
v_{i'},\quad C\in\C^\times$,
\item $x(\pi_2^+)$ satisfies IC and DC for $W(\Lambda_{i'})$,
\item $x(\pi_1)$ is maximal for $w(\mu)$, i.e. all the monomials  $x(\pi')$ that satisfy IC and DC for $L(\Lambda_i)$ and such that $w(\mu)x(\pi')v_i\neq 0$, have their $(-1)$-part $x(\pi_1')$
smaller or equal to $x(\pi_1)$.
\end{itemize}
\end{prop}

\begin{dokaz} Assume $i=0$; $\Lambda_i=\Lambda_0$, and $v_0=1=e^0$ is the highest weight vector of $L(\Lambda_0)$. Let
$$x(\pi_1)=x_{r_t s_t}(-1)\cdots x_{r_{2} s_{2}}(-1)
x_{r_1 s_1}(-1),$$ where $1\leq r_1<r_2<\dots <r_t \leq m$, $\ell
\geq s_1>s_2>\dots>s_t \geq m$. Then colors of elements of
degree $-2$ lie either below the $r_t$-th row, or left of the
$s_t$-th column (see figure \ref{DCSeq_fig}). Suppose that they lie below the $r_t$-th row. Since
$\langle \gamma_{r_p s_p},\gamma_{r_q s_q}\rangle = 0$ for $1\leq
p<q\leq t$, by \eqref{djel_vop_jed} one has
$$x(\pi_1)v_0=x_{r_t s_t}(-1)\cdots  x_{r_1
s_1}(-1) 1=C_1\cdot e^{\gamma_{r_1 s_1}+\dots+ \gamma_{r_t s_t}},$$
for some $C_1\in\C^\times$. By lemma \ref{gama_lambda}, we have
$$ x(\pi_1)v_0=C_1 \cdot e^{\lambda_{r_1}+\dots +\lambda_{r_t}+\lambda_{s_t}'+\dots +\lambda_{s_1}'}.$$

Set
$$\mu=\sum_{\begin{subarray}{c} 1\leq r<r_t \\ r \notin
\{r_1,\dots,r_t\}\end{subarray}}
\lambda_r+\sum_{\begin{subarray}{c} \ell \geq s>s_t \\ s \notin
\{s_1,\dots,s_t\}\end{subarray}}
\lambda_s'+\sum_{s=m}^{s_t-1}\lambda_s'.
$$
Weight $\mu$  is the sum of all $\lambda_r$'s, $1\leq r<r_t$, and all
$\lambda_s'$'s, $\ell\geq s\geq m$, such that in the appropriate rows and
columns, respectively, there doesn't lie any color of elements of
$x(\pi_1)$. Let $w(\mu)$ be a coefficient of $z^0=z^{\langle \mu,0
\rangle}$ in ${\mathcal Y} (e^\mu,z)$. For $\gamma\in \Gamma $,
$w(\mu)e^\gamma \neq 0$ if and only if $\langle \mu,\gamma \rangle
=0$, by \eqref{djel_vop_jed}. Because of \eqref{ScalProd_jed}, for a
monomial $x(\pi_1')$ consisting of elements of degree $-1$ and
satisfying difference conditions for $L(\Lambda_0)$, vector
$x(\pi_1')v_0$ won't be annihilated by $w(\mu)$ if and only if its
colors lie in the intersection of rows $\{r_1,\dots,r_t\}\cup
\{r_t+1,\dots,m\}$ and columns $\{s_1,\dots,s_t\}$. Clearly,
$x(\pi_1)$ is maximal among such, so if $w(\mu)x(\pi_1')v_0\neq 0$
then $x(\pi_1') \leq x(\pi_1)$.

Note that
$$
\mu+ \lambda_{r_1}+\dots +\lambda_{r_t}+\lambda_{s_t}'+\dots
+\lambda_{s_1}' = \sum_{r=1}^{r_t} \lambda_r+ \sum_{s=m}^\ell
\lambda_s' =\omega_{r_t}+\omega.
$$
Hence
\begin{eqnarray*}
w(\mu) x(\pi_1)v_0 & = & C_2 e^{\omega_{r_t}+\omega}\\
 & = & C e(\omega) v_{r_t},
\end{eqnarray*}
for some $C_2,C\in \C^\times$. Since  colors of elements of
degree $-2$ lie below the $r_t$-th row, the monomial $x(\pi_2^+)$
satisfies difference and initial conditions for $W(\Lambda_{r_t})$.
Hence the operator $w(\mu): L(\Lambda_0)\to L(\Lambda_{r_t})$
satisfies the statement of the proposition.

If colors of elements of $x(\pi)$ of degree $-2$  lie on the
left of the $s_t$-th row instead of lying below the $r_t$-th row, 
then, when constructing $\mu$, one
will replace $\lambda_s'$'s, for $m\leq s < s_t$, with $\lambda_r$'s,
for $s_t<s\leq m$. That way, we get an operator $w(\mu):
L(\Lambda_0)\to L(\Lambda_{s_t})$.

Finally, assume $1\leq i\leq \ell$; $v_i=e^{\omega_i}$ is the highest
weight vector of $L(\Lambda_i)$. Colors of elements of
$x(\pi_1)$ lie either below the $i$-th row, or on the left of $i$-th
column (see figure \ref{ICSeq_fig}). Then one constructs $\mu\in P$ similarly as before, with an
exception that if $i\leq m$, one won't take $\lambda_r$'s for
$r\leq i$, and if $i\geq m$ one won't take $\lambda_s'$'s for
$s\geq i$. For instance, if $i\leq m$ and colors of elements of
degree $-2$ in $x(\pi)$ is on the left of the $s_t$-th column, we would
set
$$\mu=\sum_{\begin{subarray}{c} i<r<r_t \\ r \notin
\{r_1,\dots,r_t\}\end{subarray}} \lambda_r+\sum_{\begin{subarray}{c}
\ell \leq s>s_t \\ s \notin \{s_1,\dots,s_t\}\end{subarray}}
\lambda_s'+\sum_{r=r_t +1}^{m}\lambda_r.$$ For the operator $w(\mu)$
we take the coefficient of $z^{\langle \mu,\omega_i\rangle}$ in
${\mathcal Y} (e^\mu,z)$. Since
$\omega_i=\lambda_1+\dots+\lambda_i$, we have
$$\mu+\gamma_{r_1 s_1}+\dots+
\gamma_{r_t s_t}+\omega_i= \omega+\omega_{s_t}.$$ Hence
$$w(\mu) x(\pi_1)v_i= C e(\omega) v_{s_t},$$
as desired.
\end{dokaz}

Proposition \ref{OpIsp1_prop} enables us to prove linear
independence of the set
$$\{x(\pi)v_i \,|\, x(\pi) \textrm{ satisfies IC and DC for } W(\Lambda_i)\}.$$
We  prove this by induction on degree and on order of monomials. The
proof is carried out simultaneously for all level $1$ standard
modules by using coefficients of intertwining operators.

Assume
\begin{equation}\label{LinZav_jed}\sum c_\pi x(\pi)v_i=0,\end{equation}
where all monomials $x(\pi)$ satisfy difference and initial
conditions for $W(\Lambda_i)$ and are of degree greater or equal to some
$-n\in\Z$. Fix $x(\pi)$ in \eqref{LinZav_jed} and suppose that
$$c_{\pi'}=0 \textrm{\quad for \quad} x(\pi')<x(\pi).$$
We want to show that $c_\pi=0$.

By proposition \ref{OpIsp1_prop}, there exists an  operator $w(\mu)$
such that
\begin{itemize}
\item $w(\mu)x(\pi_1) v_i =C e(\omega) v_{i'},\quad C \in\C^\times$,
\item $x(\pi_2^+)$ satisfies IC and DC for $W(\Lambda_{i'})$,
\item $w(\mu) x(\pi')v_i=0$\quad if \quad
$x(\pi_1')>x(\pi_1)$,
\end{itemize}
where $\Lambda_{i'}$ is another fundamental weight of $\gt$.
Applying the operator $w(\mu)$ to \eqref{LinZav_jed}  gives
\begin{eqnarray*}
0 & = & w(\mu) \sum c_{\pi'}x(\pi') v_i\\
 & = & w(\mu) \sum_{\pi_1'>\pi_1}c_{\pi'}x(\pi') v_i +w(\mu) \sum_{\pi_1'<\pi_1}c_{\pi'}x(\pi')
 v_i+ w(\mu) \sum_{\pi_1'=\pi_1}c_{\pi'}x(\pi') v_i
\end{eqnarray*}
The first sum becomes $0$ after application of $w(\mu)$, while the
second sum is also equal to $0$ by the induction hypothesis. What is
left is
\begin{eqnarray*}
0 & = & w(\mu) \sum_{\pi_1'=\pi_1}c_{\pi'}x(\pi') v_i\\
  & = & \sum_{\pi_1'=\pi_1}c_{\pi'}x(\pi_2') C e(\omega) v_{i'}\\
 & = & C e(\omega) \sum_{\pi_1'=\pi_1}c_{\pi'}x(\pi_2'^+) v_{i'}
\end{eqnarray*}
Since $e(\omega)$ is injection, it follows that
$$ \sum_{\pi_1'=\pi_1}c_{\pi'}x(\pi_2'^+)
v_{i'}=0.$$ All monomials $x(\pi_2'^+)$ satisfy difference
conditions because $x(\pi')$ were such. If some of them doesn't
satisfy initial conditions for $W(\Lambda_{i'})$, then the
corresponding monomial vectors $x(\pi_2'^+) v_{i'}$ will be equal to
$0$. Certainly, $x(\pi_2^+)$ won't be among those. We've ended up
with  a relation of linear dependence on the standard module
$L(\Lambda_{i'})$ in which all monomials are of degree greater or
equal to $-n+1$. By the induction hypothesis they are linearly
independent, and, in particular, $c_\pi=0$. We have proven

\begin{tm} \label{FSbaza_tm}
Let $L(\Lambda_i)$ be a standard $\gt$-module of level $1$. Then the
set
$$\{x(\pi)v_i\,|\,x(\pi) \textrm{ satisfies IC and DC for } W(\Lambda_i)\}$$ is a basis of $W(\Lambda_i)$.
\end{tm}

\section{Bases of standard modules}

\label{BazaStand_pogl}

Knowledge of a basis of Feigin-Stoyanovsky's type subspace
$W(\Lambda)$ was used in [P1] and [P2] to obtain a basis of the
whole standard module $L(\Lambda)$.  We're following here this
approach to obtain a basis of a standard level $1$ module
$L(\Lambda_i)$, $i=0,\dots,\ell$, for any choice of $\Z$-gradation \eqref{ZGradG_jed}.

Set
$$e=\prod_{\gamma\in\Gamma}e^\gamma=e^{\sum_{\gamma\in\Gamma}\gamma}.$$
From lemma \ref{gama_lambda} and \eqref{gama_omega_jed}, we have
$$e=e^{m\sum_{j=1}^m\lambda_j+(\ell-m+1)\sum_{j=m}^\ell\lambda_j'}=e^{(\ell+1)\omega}.$$
The following proposition was proven in [P1] and [P2] (cf. [P1,
Theorem 8.2.] or [P2, Proposition 5.2.])

\begin{prop} \label{StModGenSkup_prop} Let $L(\Lambda_i)_\mu$ be a weight subspace of
$L(\Lambda_i)$. Then there exists an integer $n_0$ such that for any
fixed $n\leq n_0$ the set of vectors
$$e^n x_{\gamma_1}(j_1)\cdots x_{\gamma_s}(j_s)v_i \in L(\Lambda_i)_\mu,$$
where  $s\geq
0,\,\gamma_1,\dots,\gamma_s\in\Gamma,\,j_1,\dots,j_s\in\Z$, is a
spanning set of $L(\Lambda_i)_\mu$. In particular,
$$L(\Lambda_i)=\langle e \rangle U(\gt_1)v_i.$$
\end{prop}

\begin{tm}
\label{StModBase_tm}
Let $L(\Lambda_i)_\mu$ be a weight subspace of a standard level $1$
$\gt$-module $L(\Lambda_i)$. Then there exists $n_0\in\Z$ such that
for any fixed $n\leq n_0$ the set of vectors
$$e^n x(\pi)v_i\in L(\Lambda_i)_\mu, \quad x(\pi)\textrm{ satisfies IC and DC for } W(\Lambda_i),$$
is a basis of $L(\Lambda_i)_\mu$. 
Moreover, for two choices of
$n_1,n_2\leq n_0$, the corresponding two bases are connected by a
diagonal matrix.
\end{tm}

\begin{dokaz}
From proposition \ref{StModGenSkup_prop} and theorem \ref{FSbaza_tm}
it follows that the set above indeed is a basis of
$L(\Lambda_i)_\mu$. It is left to prove the second part of theorem.

In order to see this, we'll find a monomial $x(\mu)\in U(\gt_1^-)$
 and $f\in \N$ such that the following holds
\begin{enumerate}
\renewcommand{\labelenumi}
{(\roman{enumi})}
\item $e(\omega)^f v_i=C x(\mu)v_i$, for some $C\in
\N$
\item $f$ divides $\ell +1$,
\item $x(\mu)$ satisfies difference and initial conditions for $W(\Lambda_i)$,
\item if a monomial $x(\pi)$ satisfies difference  and initial conditions for $W(\Lambda_i)$, then so does
a monomial $x(\pi^{-f})x(\mu)$, where $\pi^{-f}$ is a partition
defined by $\pi^{-f}(x_\gamma(-n-f))=\pi(x_\gamma(-n)),
\gamma\in\Gamma, n\in\Z$.
\end{enumerate}
Then we'll have $$ e(\omega)^f x(\pi)v_i=x(\pi^{-f})e(\omega)^f
v_i=C x(\pi^{-f})x(\mu)v_i.$$ Since $e^{\omega} x(\pi)v_i$ and
$e(\omega) x(\pi)v_i$ are proportional, the second part of the
theorem follows.

Let $x(\mu)\in U(\gt_1^-)$ be the maximal monomial satisfying
difference and initial conditions for $W(\Lambda_i)$ such that its
factors are of degree greater or equal to $-f$; we'll determine the
exact value of $f$ later. Let
$$x(\mu)=x_{p_{r},q_{r}}(-n_{r})x_{p_{r-1},q_{r-1}}(-n_{r-1})x_{p_{2},q_{2}}(-n_{2})x_{p_{1},q_{1}}(-n_{1}),$$
where factors are decreasing from right to left. The initial
conditions imply
$$x_{p_{1},q_{1}}(-n_{1})=\left\{\begin{array}{ll}
x_{1,\ell}(-1), & \textrm{ if\quad} i=0,\\
x_{1,\ell}(-2), & \textrm{ if\quad} i=m,\\
x_{i+1,\ell}(-1), & \textrm{ if\quad} 0<i<m,\\
x_{1,i-1}(-1), & \textrm{ if\quad} m<i\leq \ell.
\end{array}     \right.$$
Difference conditions between $x_{p_{t},q_{t}}(-n_{t})$ and
$x_{p_{t-1},q_{t-1}}(-n_{t-1})$ give
$$x_{p_{t},q_{t}}(-n_{t})=\left\{\begin{array}{ll}
x_{p_{t-1}+1,q_{t-1}-1}(-n_{t-1}), & \textrm{ if\quad} 1\leq p_{t-1}<m<q_{t-1}\leq\ell,\\
x_{1,q_{t-1}-1}(-n_{t-1}-1), & \textrm{ if\quad} p_{t-1}=m<q_{t-1}\leq\ell,\\
x_{p_{t-1}+1,\ell}(-n_{t-1}-1), & \textrm{ if\quad} 1\leq p_{t-1}<m=q_{t-1},\\
x_{1,\ell}(-n_{t-1}-2), & \textrm{ if\quad} p_{t-1}=m=q_{t-1}.
\end{array}     \right.$$
for $1<t\leq r$.

Degrees of elements of $x(\mu)$ are $-1,-2,\dots,-f$, respectively
from right to left. Of course, it is possible that some successive
elements are of the same degree, or that elements of a certain
degree do not occur; according to the initial and difference
conditions.

From the above observation we also see that row-indices of colors of
elements are moving cyclicly over the set
$$(1,2,\dots,m),$$ and column-indices are moving cyclicly over the
set
$$(\ell,\ell-1,\dots,m).$$
We'll choose $f$ so that we stop when we make a ``full circle'' over both sets of indices. More precisely, we choose $f$ so that the last element  $x_{p_r,q_r}(-n_r)$ of
$x(\mu)$ is
\begin{equation}
\label{PerTailF_jed}
\begin{array}{ll}
x_{m,m}(-f+1), & \textrm{ if}\quad i=0, \\
x_{m,m}(-f), & \textrm{ if}\quad i=m,\\
x_{m,i}(-f), & \textrm{ if}\quad i>m, \\
x_{i,m}(-f), & \textrm{ if}\quad 0<i<m.
\end{array}     \end{equation}
Then $r$ is equal to the
smallest common multiple of $m$ and $\ell-m+1$. From \eqref{PerTailF_jed} and proposition \ref{DCTrans_prop}, it is clear that a monomial $x(\pi)$ satisfies difference and
initial conditions for $W(\Lambda_i)$ if and only if
$x(\pi^{-f})x(\mu)$ satisfies them.

Denote by $x(\mu_j)$ the $(-j)$-part of $x(\mu)$ if there are
elements of degree $-j$ in $x(\mu)$, put $x(\mu_j)=1$ otherwise.
Suppose that $x(\mu_j)\neq 1$. Let $\gamma_j$ be the color of the
smallest element of $x(\mu_j)$. Then at least one of the indices of
$\gamma_j$ is equal to $m$. Denote by $i_j$ the other index of
$\gamma_j$ (of course, $i_j=m$ if $\gamma_j=\gamma_{m,m}$). In case
$x(\mu_j)=1$ set $i_j=0$. From \eqref{PerTailF_jed} it is obvious that $i_f=i$. The same
calculation as in the proof of proposition \ref{OpIsp1_prop} shows
that
$$x(\mu_1)v_i=C_1 e(\omega)v_{i_1},$$
and
$$x(\mu_j^{+j-1})v_{i_{j-1}}=C_j e(\omega)v_{i_j},$$
for some $C_1,\dots,C_f\in\C^\times$. Hence
\begin{eqnarray*}
x(\mu)v_i & = &x(\mu_f)\cdots x(\mu_1)v_i\\ & = & x(\mu_f)\cdots x(\mu_2)C_1 e(\omega)v_{i_1}\\ & = & C_1'e(\omega)x(\mu_f^+)\cdots x(\mu_2^+)v_{i_1}=\dots\\ & =& C e(\omega)^f v_i,\end{eqnarray*}
for some $C,C_1,C_1'\in\C^\times$.

It remains to determine $f$. If $x(\mu_j)\neq 1$ then $x(\mu_j)$
contains exactly one element whose one of the indices is equal to
$m$. If $x(\mu_j)= 1$ then $x(\mu_{j-1})$ contains $x_{m,m}(-j+1)$.
Hence number $f$ counts how many times we've crossed over $m$ while
cyclicly moving over the sets of indices $(1,2,\dots,m)$ and
$(\ell,\ell-1,\dots,m)$, i.e. $f$ is equal to the total number of
cycles we've made (over both sets of indices). Hence
$$f=\frac{r}{m}+\frac{r}{\ell-m+1}=r\frac{\ell+1}{m(\ell-m+1)}=\frac{\ell+1}{r'},$$
where $r'=\frac{m(\ell-m+1)}{r}\in\N$. In particular, $f$ divides
$\ell+1$.
\end{dokaz}

As an illustration, we can take a closer look at the case $m=1$. This is the case that was studied in [P1], for arbitrary level, and combinatorial conditions obtained there are the same as the ones that we've got.
Here, $\Gamma$ is a rectangle with $1$ row and $\ell$ columns, consisting of elements $\gamma_{11},\dots,\gamma_{1 \ell}$. Fix a fundamental weight $\Lambda_i$. A monomial $x(\pi)$ satisfies initial conditions on $W(\Lambda_i)$ if it doesn't contain elements $x_{1i}(-1),\dots,x_{1\ell}(-1)$. If we assume that elements of $x(\pi)$ are decreasing from right to left, then we can say that $x(\pi)$ satisfies difference conditions on $W(\Lambda_i)$ if for any two successive factors $x_{1s}(-j) x_{1s'}(-j')$ of $x(\pi)$ we either have $j\geq j'+2$, or $j=j'+1$ and $s<s'$. If we would write these conditions in terms of exponentials $\pi(x_\gamma(-j)), \gamma\in\Gamma, j\in \N$, we would obtain a special case of $(k,\ell+1)$-admissible  configurations, for k=1 (cf. [FJLMM],  [T]). 

We construct a periodic tail $x(\mu)$ as in the proof of theorem
\ref{StModBase_tm}. We obtain 
 $$x(\mu)=\left\{\begin{array}{ll}
x_{1,1}(-\ell)\cdots x_{1,\ell-1}(-2)x_{1,\ell}(-1), & \textrm{ if\quad} i=0,\\
x_{1,1}(-\ell-1)\cdots x_{1,\ell-1}(-3)x_{1,\ell}(-2), & \textrm{ if\quad} i=1,\\
x_{1,i}(-\ell-1)\cdots x_{1,\ell}(-i-1)x_{1,1}(-i+1)\cdots x_{1,i-1}(-1), & \textrm{ if\quad} 2\leq i\leq \ell.
\end{array}     \right.$$
which is the maximal monomial that satisfies initial and difference conditions on $W(\Lambda_i)$ and has elements of degree greater or equal to $-\ell-1$.
Since
\begin{eqnarray*}
x_{1,\ell}(-1)v_0 & = & C_0 e(\omega)v_\ell,\\
x_{1,j-1}(-1)v_j & = & C_j e(\omega) v_{j-1}, \quad j=2,\dots,\ell,\\
v_1 & = & C_1 e(\omega) v_0,
\end{eqnarray*}
for some $C_0,\dots,C_\ell\in\C^\times$, we see that
$$x(\mu)v_i=C e(\omega)^{\ell+1} v_i,$$
for some $C\in\C^\times$. 

Also, it is clear that a monomial $x(\pi)$ satisfies initial and difference cconditions on $W(\Lambda_i)$ if and only if $x(\pi^{-\ell-1})x(\mu)$ satisfies  
initial and difference conditions on $W(\Lambda_i)$.

\end{document}